\documentclass{ws-procs9x6}

\def\la{\lambda}

\DeclareMathOperator{\re}{\Re\!\text{{\small{\it e}}}}
\DeclareMathOperator{\im}{\Im\!\text{{\small{\it m}}}}

\begin{document}

\title{BLOW UP OF A CYCLE IN LOTKA-VOLTERRA TYPE EQUATIONS WITH COMPETITION-COOPERATION TERMS AND QUASI-LINEAR SYSTEMS}

\author{E. BOUSE and D. RACHINSKII$^*$}

\address{Department of Applied Mathematics, University College Cork, Ireland\\
$^*$E-mail: d.rachinskii@ucc.ie}

\begin{abstract}
We consider systems where a cycle born via the Hopf bifurcation
blows up to infinity as a parameter ranges over a finite interval.
Two examples demonstrating this effect are presented: planar
Lotka-Volterra type systems with a competition-cooperation term and
quasi-linear higher order equations.
\end{abstract}


\bodymatter

\section{Introduction}
We consider one scenario of transformation of a cycle of a
differential equation, which we call the blow up. In this scenario,
a cycle born via the Hopf bifurcation grows to infinity as a
parameter ranges over a finite segment.  From another perspective,
in the product of the phase space and the parameter axis there is a
branch of cycles connecting the equilibrium and infinity. We first
discuss the existence of such a branch for planar differential
equations where the proof can be based on the Poincare theorem. A
Lotka-Volterra type system with a competition-cooperation term is
considered as an example. Then we discuss the existence of a branch
of cycles stretching from zero to infinity for a class of higher
order quasilinear equations: this theorem continues the results of
Refs.~\refcite{Branch,Chaos}. The results for planar systems are
presented in the next section. Section 3 contains the main result
for higher order equations. We briefly sketch some points of the
proofs.

\section{Planar systems}
Consider a planar system
\begin{equation}\label{planar}
x'=f(x,y;\lambda),\quad y'=g(x,y;\lambda)
\end{equation}
with a parameter $\lambda\in[0,1]$, where the functions $f, g$ are
continuously differentiable with respect to the phase variables
$x,y$ and continuous with respect to the set of all arguments. We
say that the system has a branch of cycles connecting the origin and
infinity if for any bounded open domain $G\ni (0,0)$ of the phase
plane the system has a cycle (for some $\lambda\in [0,1]$) that
belongs to the closure $\bar G$ of $G$ and touches the boundary
$\partial G$ of $G$. This definition, close to the classical weak
definition of continuous branches of fixed points \cite{Krasn}\,, is
discussed in Ref.~\refcite{Branch}.

\begin{proposition}\label{p1}
Suppose that the following conditions hold:

\begin{romanlist}
\item
$f(0,0;\lambda)=0$ and the origin is the only equilibrium point of
system (\ref{planar}) for all $\lambda\in[0,1]$;

\item
the Jacobi matrix $J=J(\lambda)$ of system (\ref{planar}) at the
origin is invertible for all $\lambda\in[0,1]$;

\item
the spectrum of $J$ is in the open right half-plane for $\lambda=0$
and in the open left half-plane for $\lambda=1$;

\item
the system does not have cycles for $\lambda=0$ and $\lambda=1$.
\end{romanlist}
Then system \eqref{planar} has a branch of cycles connecting the
origin and infinity.
\end{proposition}

Condition (ii) ensures that no equilibrium branches from the origin,
which agrees with condition (i). According to condition (iii), the
origin changes stability as $\lambda$ ranges over the segment
$[0,1]$. Conditions (ii) and (iii) guarantee that the
Andronov\,--\,Hopf bifurcation occurs within this segment. Condition
(iv) is satisfied, for example, if the system has a global Lyapunov
function for $\lambda=0$ and $\lambda=1$. Under the assumptions of
Proposition \ref{p1} the system does not have homoclinic orbits; the
cycles surround the origin.

As an example, consider the Lotka-Volterra type equations
\begin{equation}\label{LV}
x'=x(a-by),\quad y'=y(-c+dx+f(y;\lambda));\qquad x, y>0
\end{equation}
with positive parameters $a,b,c,d$. The last term $f$ in the
predator equation accounts for competition ($f<0$) or cooperation
($f>0$) in the predator population. We assume that the behavior of
the predator can depend on its number: it competes if the population
number is above a certain threshold $\bar y=\bar y(\lambda)$ and
starts to cooperate when the population falls below the threshold,
hence $f(y;\lambda)<0$ for $y> \bar y(\lambda)$ and $f(y;\lambda)>0$
for $y< \bar y(\lambda)$. For example, $f=\arctan y-\lambda y$ with
$\lambda\in(0,1)$.

Suppose that $c>f(ab^{-1};\lambda)$ for all $\lambda$ (the condition
$c>f$ means that the predator extincts in the absence of prey). Then
system (\ref{LV}) has a unique positive equilibrium
$(x_{*},y_{*})=(d^{-1}(c-f(ab^{-1};\lambda)), ab^{-1})$ and
Proposition \ref{p1} can be applied after the logarithmic coordinate
transformation and the shift of the equilibrium to zero. The
equilibrium $(x_*,y_*)$ is  stable if $f'_y <0$ and unstable if $
f'_y
>0$ at $y=y_{*}$, where $f'_y=\partial f/\partial y$ is the partial derivative of $f=f(y;\lambda)$.
Hence conditions (i) -- (iii) of Proposition \ref{p1} are satisfied
if
\begin{equation}\label{iiii}
f'_y (y_{*};0)>0, \qquad f'_y (y_{*};1)<0.
\end{equation}
If in addition
\begin{eqnarray}\label{iv}
(y-y_{*})(f(y;0)-f(y_{*};0))>0,& y>0,\ y\ne y_{*},\\
(y-y_{*})(f(y;1)-f(y_{*};1))<0,& y>0,\ y\ne y_{*},\label{iv'}
\end{eqnarray}
then condition (iv) is also satisfied, because in this case $
V=(x-x_{*}\ln x)d+(y-y_{*}\ln y)b $ is a Lyapunov function of system
(\ref{LV}) for $\lambda=0,1$ with $\dot
V=(y-y_{*})(f(y;\lambda)-f(y_{*};\lambda))b$. Hence, relations
(\ref{iiii}) -- (\ref{iv'}) ensure the existence of a branch of
cycles connecting the equilibrium and infinity for system (\ref{LV})
(the cycles lie in the positive quadrant $x,y>0$ where the system is
defined). In particular, these relations hold if $f$ strictly
increases for $\lambda=0$ and strictly decreases for $\lambda=1$, as
in the above example $f=\arctan y-\lambda y$. Numerical simulations
confirm that a stable positive cycle born via the Hopf bifurcation
blows up to infinity for this $f$ and demonstrate the same effect
for other competition-cooperation terms, such as $f=y-\lambda y^{2}$
or $f=y^{2}-\lambda y^{3}$, included in the equation for predator,
or prey, or both.

To prove Proposition \ref{p1}, one can first note that under its
conditions the equilibrium can not have eigenvalues of different
sign and consequently the system does not have homoclinic orbits.
Using the Poincare theorem, one derives from this fact that if the
system has an orbit $\gamma$ in a bounded domain $\bar G$ with
$\gamma\cap\partial G\ne\emptyset$, then it also has a cycle
$\mathcal{C}\subset \bar G$ with $\mathcal{C}\cap\partial
G\ne\emptyset$. To complete the proof by contradiction, assume that
there is no such a cycle and hence no such an orbit $\gamma$ for
some bounded domain $G\ni 0$. Consequently, $S_\lambda\cap \partial
G=\emptyset$, where $S_\lambda$ denotes the invariant set of the
system in the domain $\bar G$ (note that $0\in S_\lambda$).
Therefore $\bar G$ is an isolating neighborhood for $S_\lambda$ and
the Conley index ${\rm Ind}\, S_\lambda$ of $S_\lambda$ with respect
to $\bar G$ is defined \cite{Conley}\,. Because the system has no
homoclinic orbits and, by condition (iv), there is no cycles for
$\lambda=0,1$, the Poincare theorem implies that $S_{0}=S_{1}=0$.
Moreover, condition (iii) implies ${\rm Ind}\, S_0\ne {\rm Ind}\,
S_1$. This, however, contradicts the invariance of the Conley index
under homotopic transformation of the vector field: ${\rm Ind}\,
S_\lambda$ should be the same for all $\lambda$ for any isolating
neighbprhood $\bar G$ of $S_{\lambda}$. Given any open bounded $G\ni
0$, this contradiction proves the existence of a cycle
$\mathcal{C}\subset \bar G$ with $\mathcal{C}\cap\partial
G\ne\emptyset$ for some $\lambda$, i.e. the conclusion of the
proposition.

\section{Quasi-linear higher order equations}

Consider the equation
\begin{equation}\label{1}
L\left(\frac{d}{dt};\lambda\right)x=f(x;\la),
\end{equation}
where $L(p;\lambda)=p^\ell
+a_{1}(\la)p^{\ell-1}+\cdots+a_0(\lambda)$ is a polynomial with
continuously differentiable coefficients. Assume that the continuous
function $f(x;\la)$ satisfies $f(0;\la)\equiv 0$ and the global
Lipschitz estimates
\begin{equation}\label{lipx}
%$$
|f(x_{1};\la)-f(x_{2};\la)| \le k|x_{1}-x_{2}|,\
%\end{equation}
%\begin{equation}\label{lipl}
|f(x;\la_{1})-f(x;\la_{2})| \le l|x||\la_{1}-\la_{2}|;
%$$
\end{equation}
%A simple example of such function is $f(x,\la)=(2+\sin \la)\,x\sin\ln |x|$.
hence the equation has the zero solution $x\equiv 0$ for all
$\lambda$. Define the matrix
%\begin{equation}\label{J}
$$
J(w,\lambda)=\left(
\begin{array}{rrrrr}
\re L_{\la}'(\la;wi)&&&& -\im L_{p}'(\la;wi)\\
\im L_{\la}'(\la;wi)&&&& \re L_{p}'(\la;wi)
\end{array}
\right).
$$
%\end{equation}
%where $L_{\la}'$, $L_{p}'$ are the partial derivatives of $L=L(p;\lambda)$.
%Denote by $\Z$  the set of all nonnegative integers, $\Z=\{0,1,2,\ldots\}$.
We say that equation \eqref{1} has a Lipschitz continuous branch of
cycles connecting zero and infinity if there are Lipschitz
continuous functions $\lambda(r)$, $w(r)$ with values in segments
$[\la_-,\la_+]$, $[w_-,w_+]$ ($w_->0$) such that for every $r>0$
equation (\ref{1}) with $\lambda=\lambda(r)$ has a periodic solution
$x_r(t)=x(t;r)$ of the period $2\pi/w(r)$, the function $x(t;r)$ is
Lipschitz continuous in $r$ and
$$
\|x_r\|_C\to 0\ \ \text{as}\ \ r\to 0,\qquad \|x_r\|_C\to \infty\ \
\text{as}\ \ r\to\infty.
$$

\begin{theorem}\label{t1}
Assume that for some $q>0$ the relation $|L(wi;\la)|\le q$ defines a
simply connected bounded domain $D_{q}$ on the plane $(w,\la)$, the
equation $L(wi;\la)=0$ has a unique solution $(w_0,\la_{0})$ in
$D_{q}$, the matrix $J(w,\lambda)$ is nondegenerate in $D_{q}$, and
$L( nwi;\lambda)\ne 0$ in $D_q$ for any integer $ n\ne\pm 1$. Then
there are sufficiently small $k, l>0$ such that equation \eqref{1}
with any function $f$ satisfying the estimates \eqref{lipx} has a
Lipschitz continuous branch of cycles connecting zero and infinity.
\end{theorem}

The method of the proof of Theorem \ref{t1} leads to explicit
estimates of the Lipschitz coefficients $k,l$, which ensure the
existence of the branch of cycles connecting zero and infinity.

A natural parameter $r$ is the amplitude of the first harmonics of
the periodic solution. Theorem \ref{t1} can be proved by contraction
mapping principle. To construct the corresponding mapping, let us
first note that for any $(w,\lambda)\in D_q$ the differential
operator $L(w \frac{d}{dt};\lambda)$ with the $2\pi$-periodic
boundary conditions is invertible on the codimension 2 subspace
$\mathbb{E}$ of $\mathbb{L}^2=\mathbb{L}^2(0,2\pi)$ which is
orthogonal to $\sin t$ and $\cos t$ (this operator, however, is not
invertible on the whole space $\mathbb{L}^2$, because
$L(w_0i;\lambda_0)=0)$. Secondly, the planar map $(w,\lambda)
\mapsto (\re L(wi;\lambda),\im L(wi;\lambda))$ is invertible on
$D_q$. Now, denote by $P_s$ and $P_c$ the orthogonal projectors onto
$\sin t$ and $\cos t$ in $\mathbb{L}^2$, define the orthogonal
projector $Q=I-P_s-P_c$ onto $\mathbb{E}=Q\mathbb{L}^2$ and consider
the space of triples $(u,v,y)\in \mathbb{R}\times\mathbb{R}\times
Q\mathbb{L}^2$ with the norm
$\|(u,v,y)\|=\sqrt{u^2+v^2+\|y\|^2_{L^2}}$, where $y=y(t)$. In this
space, for each value of the parameter $r>0$ consider the mapping
$$
(u,v,y)\mapsto A_r(u,v,y)=r^{-1}(P_s f(x(t);\lambda), P_c
f(x(t);\lambda), Qf(x(t);\lambda)),
$$
where $x=x(t)$ and $\lambda$ are defined by the relations
\begin{eqnarray}\label{.}
u=\re L(wi;\lambda),\ \ v=\im L(wi;\lambda),\ \
x(t)=r(\pi^{-1/2}\sin t +h(t)), \\
\label{..}
 L\Bigl(w\frac{d}{dt};\lambda\Bigr) h(t)=y(t),\ \ h(0)=h(2\pi),
 h'(0)=h'(2\pi),\  h\in Q{\mathbb L}^2.
\end{eqnarray}
Due to the invertibility of these relations mentioned above, the
mapping $A_r$ is well-defined for all $u^2+v^2\le q^2$ and all $y\in
Q\mathbb{L}^2$. The definition of $A_r$  ensures that every fixed
point of $A_r$ defines a $2\pi$-periodic solution $x=x(t)$ of the
equation $L(w\frac{d}{dt};\lambda)x=f(x;\lambda)$ by the formulas
\eqref{.}, \eqref{..} and hence a $2\pi/w$-periodic solution $x(wt)$
of equation \eqref{1}. The proof is completed by showing that if $k,
l$ are sufficiently small then $A_r$ is a contraction on the ball
$\|(u,v,y)\|\le q$, which is invariant for $A_r$, for each $r>0$.

%Set
%\begin{equation}\label{In}
%I_{n}=\left(\begin{array}{cc} 1&0\\ 0&n
%\end{array}\right).
%\end{equation}
%We shall use the norm
%\begin{equation}\label{norm}
%\|G\|=\max_{0\le\theta\le2\pi}\sqrt{(g_{11}\cos\theta+g_{12}\sin\theta)^{2}+
%(g_{21}\cos\theta+g_{22}\sin\theta)^{2}}
%\end{equation}
%of the real matrixes
%$$
%G=\left(\begin{array}{cc} g_{11}&g_{12}\\ g_{21}&g_{22}
%\end{array}\right).
%$$

%=================================================================================

%Suppose that
%\begin{equation}\label{cond}
%\begin{array}{ccl}
%1&>& \displaystyle\max_{|L(\la,wi)|\le q;\ n\in\Z,\,n\ne1}
%\frac{k^{2}}{|L(\la,nwi)|^{2}}
%\\
%\\
%&+ &\displaystyle \left(\max_{|L(\la,wi)|\le q;\ n\in\Z,\,n\ne1}
%\frac{k\|J(\la,nw)I_{n}J^{-1}(\la,w)\|}{|L(\la,nwi)|^{2}}\sqrt{
%q^{2}-|L(\la,wi)|^{2}}\right.
%\\
%\\
%&+ &\displaystyle\left.\max_{|L(\la,wi)|\le q;\ n\in\Z,\,n\ne1}
%\frac{l|L'_{p}(\la,wi)|}{|\det J(\la,w)|} \sqrt{ 1 + \frac{{
%q^{2}-|L(\la,wi)|^{2}}}{|L(\la,nwi)|^{2}}}\ \right)^{2}
%\end{array}
%\end{equation}
%and
%\begin{equation}\label{cond0}
%\max_{|L(\la,wi)|\le q;\ n\in\Z,\, n\ne1} k\sqrt{1
%+\frac{q^{2}-|L(\la,wi)|^{2}}{|L(\la,nwi)|^{2}}}\le q.
%\end{equation}

\section*{Acknowledgements}
This publication has emanated from research conducted with the
financial support of Science Foundation Ireland and IRCSET.

\end{document}